\documentclass[10pt]{amsart}

\usepackage[english]{babel}
\usepackage[T1]{fontenc}
\usepackage[utf8]{inputenc}
\usepackage{amsmath}
\usepackage{amssymb}
\usepackage{amsthm}
\usepackage{color}
\usepackage[a4paper]{geometry}
\usepackage{tikz}
\usepackage{mathdots}
\usepackage[all]{xy}
\usepackage{hyperref}
\usepackage{accents}
\usepackage{lmodern}
\usepackage[multiple]{footmisc}
\usepackage{mathrsfs}

\theoremstyle{definition}
\newtheorem{definition}{Definition}[section]

\newtheorem{conj}[definition]{Conjecture}

\theoremstyle{plain}
\newtheorem{theorem}[definition]{Theorem}
\newtheorem{prop}[definition]{Proposition}
\newtheorem{cor}[definition]{Corollary}
\newtheorem{lem}[definition]{Lemma}

\newcommand{\mb}{\mathbb}
\newcommand{\mbf}{\mathbf}
\newcommand{\mc}{\mathcal}

\newcommand{\bs}{\boldsymbol}

\begin{document}

\title{Multiplicatively badly approximable matrices up to logarithmic factors}

\author{Reynold Fregoli}
\address{Department of Mathematics\\ 
Royal Holloway, University of London\\ 
TW20 0EX Egham\\ 
UK}
\email{Reynold.Fregoli.2017@live.rhul.ac.uk}

\date{\today, and in revised form ....}

\dedicatory{}


\begin{abstract}
Let $\|x\|$ denote the distance from $x\in\mb{R}$ to the nearest integer. In this paper, we prove an existence and density statement for matrices $\bs{A}\in\mb{R}^{m\times n}$ satisfying
$$\liminf_{|\bs{q}|_{\infty}\to +\infty}\prod_{j=1}^{n}\max\{1,|q_{j}|\}\log\left(\prod_{j=1}^{n}\max\{1,|q_{j}|\}\right)^{m+n-1}\prod_{i=1}^{m}\|A_{i}\bs{q}\|>0,$$
where the vector $\bs{q}$ ranges in $\mb{Z}^{n}$ and $A_{i}$ are the rows of the matrix $\bs{A}$. This result extends a previous result of Moshchevitin for $2$-dimensional vectors to arbitrary dimension. The estimates needed to apply Moshchevitin's method to the case $m>2$ are not currently available. We therefore develop a substantially different method, that allows us to overcome this issue. We also generalise this existence result to the inhomogeneous setting. Matrices with the above property appear to have a very small sum of reciprocals of fractional parts. This fact helps us to shed light on a question raised by L\^e and Vaaler, thereby proving some new estimates for such sums in higher dimension.
\end{abstract}

\maketitle

\section{Introduction}

\subsection{Notation}
For $x\in\mb{R}$ we denote by $\|x\|$ the distance from $x$ to the nearest integer. 
For a matrix $\bs{A}\in\mb{R}^{m\times n}$ we denote by $A_{i}\in\mb{R}^{n}$ the rows of $\bs{A}$ ($i=1,\dotsc,m$). If $\bs{x}\in\mb{R}^{n}$, we denote by $A_{i}\bs{x}$ the sum $\sum_{j=1}^{n}A_{ij}x_{j}$. Given a set $X$ and a pair of functions $f,g:X\to\mb{R}$, we write $f\ll g$ (or $f\gg g$) if there exists a constant $c>0$ such that $f(x)\leq cg(x)$ (or $f(x)\geq cg(x)$) for all $x\in X$. If the constant $c$ depends on some parameters, we write them under the symbol $\ll$ (or $\gg$). We denote by $|\cdot|_{2}$ the Euclidean norm and by $|\cdot|_{\infty}$ the supremum norm in $\mb{R}^{n}$. We also denote by $\textup{dist}_{2}$ and $\textup{dist}_{\infty}$ the Euclidean and supremum distances respectively. For a set $X\subset\mb{R}^{n}$ we denote by $\textup{diam}(X)$ its diameter and by $\textup{Vol}(X)$ its $n$-dimensional Hausdorff measure. If $X$ is a (hyper)cube, we denote by $\textup{edge}(X)$ the length of its edges, i.e., its $1$-dimensional faces. If $f:\mb{Z}^{n}\to[0,+\infty)$ is a function, we denote by $\liminf_{|\bs{q}|_{\infty}\to +\infty}f(\bs{q})$ the number $\liminf_{q\to +\infty} \min\{f(\bs{q}):|\bs{q}|_{\infty}=q\}$. Finally, we intend any product $\prod_{a}^{b}$ where $b<a$ as $1$.

\subsection{Background}
It is well known that the set of real numbers $\alpha\in\mb{R}$ such that
\begin{equation}
\label{eq:badappnumbers}
\liminf_{q\to\infty}q\|q\alpha\|>0
\end{equation}
is non empty and has full Hausdorff dimension. Such numbers are called badly approximable. The notion of bad approximablity can be extended to a higher-dimensional setting. In particular, a matrix $\bs{A}\in\mb{R}^{m\times n}$ (where $m,n$ are positive integers) is said to be badly approximable if, in analogy to (\ref{eq:badappnumbers}),
\begin{equation}
\liminf_{|\bs{q}|_{\infty}\to +\infty}|\bs{q}|_{\infty}^{n}\max_{i=1}^{m}\{\|A_{i}\bs{q}\|\}^{m}>0.\nonumber
\end{equation}
Schmidt \cite{Schmidt:BadlyApprox} showed that the set of such matrices has full Hausdorff dimension in $\mb{R}^{m\times n}$.

A multiplicative generalisation of (\ref{eq:badappnumbers}) has also extensively been studied \cite{Bugeaud:Multiplicative}. A matrix $\bs{A}\in\mb{R}^{m\times n}$ is said to be multiplicatively badly approximable if
\begin{equation}
\label{eq:multbad}
\liminf_{|\bs{q}|_{\infty}\to +\infty}\prod_{j=1}^{n}\max\{1,|q_{j}|\}\prod_{i=1}^{m}\|A_{i}\bs{q}\|>0.
\end{equation}
To simplify the notation, throughout this paper we write
$$\prod(\bs{q}):=\prod_{j=1}^{n}\max\{1,|q_{j}|\}$$
for all $\bs{q}\in\mb{Z}^{n}$.

Proving the existence of multiplicatively badly approximable matrices is a major problem in Diophantine approximation. The famous Littlewood conjecture states that for any pair of real numbers $\alpha,\beta\in\mb{R}$ it holds
\begin{equation}
\liminf_{q\to\infty}q\|q\alpha\|\|q\beta\|=0,
\end{equation}
or that, in other words, there are no $2\times 1$ multiplicatively badly approximable matrices. If the Littlewood conjecture were true, there would be no multiplicatively badly approximable matrices (except for $n=m=1$). This follows from the fact that every submatrix of a multiplicatively badly approximable matrix is itself multiplicatively badly approximable and from the transference principle (see \cite[Theorem 2.2]{LeVaaler:Sumsof}).

Since the set of multiplicatively badly approximable matrices with this definition could be empty, some authors (Badziahin, Velani, etc.) started working on a different definition of multiplicative bad approximability. Their idea was to weaken the Diophantine condition in (\ref{eq:multbad}) to allow for more flexibility. One possible way of doing this could be increasing the exponent of the factor $\prod(\bs{q})$. This modification however, adds too many matrices to the set in consequence of the two following $0-1$ results. 

\begin{theorem}[Gallagher]
\label{thm:Gallagher}
Let $m$ be a positive integer and let $\psi:\mb{N}\to(0,1]$ be a non-increasing\footnote{We say that a function $f:A\to B$ with $A,B\subset\mb{R}$ is non-increasing if $f(x)\geq f(y)$ for all $x<y$.} function. Let also
$$W^{\times}(m,1,\psi):=\left\{\bs{A}\in[0,1]^{m\times 1}:\prod_{i=1}^{m}\|A_{i}q\|<\psi\left(|q|\right)\ \mbox{for infinitely many }q\in\mb{Z}\right\}.$$
Then,
$$
\mathscr{L}(W^{\times}(m,1,\psi))=\begin{cases}
0 & \mbox{if }\sum_{q=1}^{+\infty}\psi\left(q\right)\log\left(\psi(q)^{-1}\right)^{m-1}< +\infty \\
1 & \mbox{if }\sum_{q=1}^{+\infty}\psi\left(q\right)\log\left(\psi(q)^{-1}\right)^{m-1}= +\infty
\end{cases},
$$
where $\mathscr{L}$ stands for the $m$-dimensional Lebesgue measure.
\end{theorem}

\begin{theorem}[Sprind\v{z}uk]
\label{thm:Sprindzuk}
Let $m,n$ be positive integers and let $\psi:\mb{N}\to(0,1]$ be any function. Let also
$$W^{\times}(m,n,\psi):=\left\{\bs{A}\in[0,1]^{m\times n}:\prod_{i=1}^{m}\|A_{i}\bs{q}\|<\psi\left(\prod(\bs{q})\right)\ \mbox{for infinitely many }\bs{q}\in\mb{Z}^{n}\right\}.$$
Then,
$$
\mathscr{L}(W^{\times}(m,n,\psi))=\begin{cases}
0 & \mbox{if }\sum_{\bs{q}\in \mb{Z}^{n}}\psi\left(\prod(\bs{q})\right)\log\left(\psi\left(\prod(\bs{q})\right)^{-1}\right)^{m-1}< +\infty \\
1 & \mbox{if }\sum_{\bs{q}\in S}\psi\left(\prod(\bs{q})\right)\log\left(\psi\left(\prod(\bs{q})\right)^{-1}\right)^{m-1}= +\infty
\end{cases},
$$
where $\mathscr{L}$ stands for the $mn$-dimensional Lebesgue measure, and $S$ is any infinite set of pairwise linearly independent vectors in $\mb{Z}^{n}$.
\end{theorem}

Theorem \ref{thm:Gallagher} follows from the more general \cite[Theorem 1]{Gallagher:MetricSimultaneous}, whereas Theorem \ref{thm:Sprindzuk} follows from its multi-dimensional analogue \cite[Chapter 1, Theorem 13]{Sprindzuk:Metric}. Note that there is a discrepancy between the cases $n=1$ and $n>1$. In particular, Theorem \ref{thm:Sprindzuk} does not imply Theorem \ref{thm:Gallagher}, since for $n=1$ there are no infinite subsets of pairwise linearly independent vectors in $\mb{Z}$. 

Gallagher and Sprind\v{z}uk's Theorems both imply that the set of matrices $\bs{A}\in\mb{R}^{m\times n}$ such that
$$\liminf_{|\bs{q}|_{\infty}\to +\infty}\prod(\bs{q})^{1+\varepsilon}\prod_{i=1}^{m}\|A_{i}\bs{q}\|>0,$$
has full Lebesgue measure in $\mb{R}^{m\times n}$ for all $\varepsilon>0$. Therefore, a finer indicator (in comparison with the exponent of $\prod(\bs{q})$) is required to find a non-empty set that would not coincide with almost all the space $\mb{R}^{m\times n}$. A natural approach is to allow for logarithmic factors, i.e., to consider the set
\begin{equation}
\textup{Mad}^{\lambda}(m,n):=\left\{\bs{A}\in\mb{R}^{m\times n}:\liminf_{|\bs{q}|_{\infty}\to +\infty}\prod(\bs{q})\log\left(\prod(\bs{q})\right)^{\lambda}\prod_{i=1}^{m}\|A_{i}\bs{q}\|>0\right\}.
\end{equation} 

It follows from Theorems \ref{thm:Gallagher} and \ref{thm:Sprindzuk} that $\textup{Mad}^{\lambda}(m,n)$ has full Lebesgue measure for $\lambda>m+n-1$ and zero Lebesgue measure for $\lambda\leq m+n-1$. However, it could happen, for example, that $\textup{Mad}^{\lambda}(m,n)$ is empty for $\lambda\leq m+n-1$, and this is the case that we treat in this paper. 

To have a better understanding of the situation, we look at the analogue of the set $\textup{Mad}^{\lambda}(m,n)$ in the standard setting, i.e., the set
\begin{equation}
\textup{Bad}^{\lambda}(m,n):=\left\{\bs{A}\in\mb{R}^{m\times n}:\liminf_{|\bs{q}|_{\infty}\to +\infty}|\bs{q}|_{\infty}^{n}\log(|\bs{q}|_{\infty})^{\lambda}\max\{\|A_{1}\bs{q}\|,\dotsc,\|A_{m}\bs{q}\|\}^{m}>0\right\}.
\end{equation} 

The analogue of Theorems \ref{thm:Gallagher} and \ref{thm:Sprindzuk} in the additive setting is the Khintchine-Groshev Theorem (see references in \cite{Beresnevich:Classic}), which we report for the convenience of the reader. 

\begin{theorem}[Khintchine-Groshev]
\label{thm:Khintchine-Groshev}
Let $m,n$ be positive integers and let $\psi:\mb{N}\to(0,1]$ be a non-increasing function. Let also
$$W^{+}(m,n,\psi):=\left\{\bs{A}\in[0,1]^{m\times n}:\max_{i=1}^{m}\{\|A_{i}\bs{q}\|\}^{m}<\psi(|\bs{q}|_{\infty})\ \mbox{for infinitely many }\bs{q}\in\mb{Z}^{n}\right\}.$$ 
Then,
$$
\mathscr{L}(W^{+}(m,n,\psi))=\begin{cases}
0 & \mbox{if }\sum_{q=1}^{+\infty}\psi(q)q^{n-1}< +\infty \\
1 & \mbox{if }\sum_{q=1}^{+\infty}\psi(q)q^{n-1}= +\infty
\end{cases},
$$
where $\mathscr{L}$ stands for the $mn$-dimensional Lebesgue measure.
\end{theorem}

This theorem, in combination with Schmidt's dimensional result for badly approximable matrices \cite{Schmidt:BadlyApprox} and Dirichlet's Theorem, implies that 
\begin{equation}
\label{eq:badmetric}
\textup{Bad}^{\lambda}(m,n)=\begin{cases}
\emptyset & \mbox{if }\lambda< 0 \\
\mbox{full Hausdorff dimension but zero Lebesgue measure set} & \mbox{if }0\leq\lambda\leq 1 \\
\mbox{full Lebesgue measure set} & \mbox{if }\lambda> 1
\end{cases}.
\end{equation}
\noindent In particular, we observe that "shaving off" a logarithm factor from the Lebesgue $0-1$ "switch over" gives the set of badly approximable matrices. 

Let us move back to the multiplicative framework and draw a comparison. To keep things simple we set $m=2$, $n=1$. We note that Theorem \ref{thm:Khintchine-Groshev} for $m=2$, $n=1$ and Theorem \ref{thm:Gallagher} for $m=2$ differ only by the presence of a logarithmic factor in the sum. In particular, Gallagher's Theorem implies that
$$\mathscr{L}\left(\textup{Mad}^{\lambda}(2,1)\right)=\begin{cases}
0 & \mbox{if }\lambda\leq 2 \\
+\infty & \mbox{if }\lambda>2
\end{cases}.
$$ 
Drawing inspiration from (\ref{eq:badmetric}) and from the "shaving off" phenomenon, Badziahin and Velani \cite[Statements L1-L3]{BadziahinVelani:MultiplicativelyBadly} made the following conjecture.
\begin{conj}[Badziahin-Velani]
\label{conj:BadzVel}
$$\textup{Mad}^{\lambda}(2,1)=
\begin{cases}
\emptyset & \mbox{if }\lambda<1 \\
\mbox{full Hausdorff dimension but zero Lebesgue measure set} & \mbox{if } 1\leq\lambda\leq 2 \\
\mbox{full Lebesgue measure set} & \mbox{if }\lambda > 2
\end{cases}.
$$
\end{conj}
This is also supported by heuristic volume arguments of Peck \cite{Peck:Simultaneous}, and Pollington and Velani \cite{PollingtonVelani:OnaProblem} (see references in \cite{BadziahinVelani:MultiplicativelyBadly}). If this conjecture were true, the set $\textup{Mad}^{1}(2,1)$ would be rightfully regarded as the multiplicative analogue of the set $\textup{Bad}^{0}(2,1)$ (i.e., the set of badly approximable vectors in $\mb{R}^{2}$). Note that Conjecture \ref{conj:BadzVel} implies the Littlewood Conjecture.

Multiple authors have contributed towards a partial solution of \ref{conj:BadzVel}. Moshchevitin \cite{Moshchevitin:BadlyApp} was the first to show that $\textup{Mad}^{2}(2,1)\neq\emptyset$ by using the so-called Peres-Schlag method. Subsequently, Bugeaud and Moschevitin \cite{BugeaudMoshchevitin:BadlyApp} showed that $\textup{dim}\,\textup{Mad}^{2}(2,1)=2$, where $\textup{dim}$ denotes the Hausdorff dimension. Finally, Badziahin \cite{Badziahin:OnMultiplicatively} showed that $\textup{dim}\,\textup{Mad}^{\lambda}(2,1)=2$ for all  $\lambda>1$. The case $\lambda=1$ of this conjecture is still unsolved.

\subsection{Main result}

In analogy with Moshchevitin's result \cite{Moshchevitin:BadlyApp}, we show in this paper that the set $\textup{Mad}^{m+n-1}(m,n)$ is dense and uncountable, in particular, it is non-empty for all $m,n\in\mb{N}$. We furthermore generalise this result to the inhomogeneous setting.

Let $C\subset \mb{R}^{m\times n}$ be a cube of edge $\ell$. For $f:[0,+\infty)\to[1,+\infty)$ non-decreasing\footnote{We say that a function $f:A\to B$ with $A,B\subset\mb{R}$ is non-decreasing if $f(x)\leq f(y)$ for all $x< y$.}, $\bs{\gamma}\in\mb{R}^{m}$, and $c>0$ we consider the set
\begin{multline}
\textup{Mad}_{m,n}(C,\bs{\gamma},f,c):=\left\{\bs{A}\in C:\prod(\bs{q})\|A_{1}\bs{q}+\gamma_{1}\|\dotsm\|A_{m}\bs{q}+\gamma_{m}\|>\frac{c}{f(\prod(\bs{q}))}\right. \\
\mbox{for all } \bs{q}\in\mb{Z}^{n}\setminus\{\bs{0}\}\bigg\}.\nonumber
\end{multline}
For $x\in[0,+\infty)$ we set $\log^{*}(x):=\log\left(\max\{e,x\}\right)$, where $e=2.71828\dots$ is the base of the natural logarithm. Then, the following holds.

\begin{prop}
\label{prop:mainprop}
Let $m,n\in\mb{N}$ with $m+n\geq 3$, let $C$ be a cube in $\mb{R}^{m\times n}$, and let $\bs{\gamma}\in\mb{R}^{m}$. Then, there exists a constant $c=c(m,n,\ell)>0$, only depending on the integers $m$ and $n$, and the edge of the cube $C$, such that for any countable (possibly finite) family of hyperplanes $\mathscr{H}$ lying in $\mb{R}^{m\times n}$ we have  
$$\textup{Mad}_{m,n}\left(C,\bs{\gamma},\log^{*}(x)^{m+n-1},c\right)\setminus \bigcup_{H\in\mathscr{H}}H\neq\emptyset.$$
\end{prop}

\noindent Proposition \ref{prop:mainprop} immediately implies the following corollary.

\begin{cor}
Let $m+n\geq 3$. Then, for all $\bs{\gamma}\in\mb{R}^{m}$ the set
$$\textup{Mad}^{m+n-1}(m,n,\bs{\gamma}):=\left\{\bs{A}\in\mb{R}^{m\times n}:\liminf_{|\bs{q}|_{\infty}\to +\infty}\prod(\bs{q})\log\left(\prod(\bs{q})\right)^{m+n-1}\prod_{i=1}^{m}\|A_{i}\bs{q}+\gamma_{i}\|>0\right\}$$
is everywhere dense in $\mb{R}^{m\times n}$ and does not lie on a countable union of hyperplanes.
\end{cor}

Note that, for certain values of $\bs{\gamma}\in\mb{R}^{m}$ the fact that the set $\textup{Mad}^{m+n-1}(m,n,\bs{\gamma})$ is uncountable is trivial (e.g., take $n=1$, $\gamma_{1},\dotsc,\gamma_{m-1}\notin\mb{Q}$, $\gamma_{m}=0$, $A_{1},\dotsc,A_{m-1}\in\mb{Z}$, and $A_{m}$ badly approximable). However, the fact that $\textup{Mad}^{m+n-1}(m,n,\bs{\gamma})$ does not lie on a countable union of hyperplanes implies that there exist matrices $\bs{A}\in\textup{Mad}^{m+n-1}(m,n,\bs{\gamma})$ whose entries, along with $1$ and the entries of $\bs{\gamma}$, are linearly independent over $\mb{Q}$. This excludes many of the most trivial examples. That said, for $\bs{\gamma}=\bs{0}$ even the non-emptiness of the set $\textup{Mad}^{m+n-1}(m,n,\bs{\gamma})$ appears to be non-trivial.

It is worth observing that to prove Proposition \ref{prop:mainprop} we do not follow the Peres-Schlag method, i.e., the method used by Moshchevitin to show that $\textup{Mad}^{2}(2,1)\neq\emptyset$ (see \cite{Moshchevitin:BadlyApp}). Moshchevitin's proof relies both on the one dimensional case ($m=n=1$), and on estimates for the sum
$$\sum_{q=1}^{Q}\frac{1}{q\|q\alpha\|}.$$ 
This sum is known to grow like $O(\log(Q)^{2})$ for almost all $\alpha\in\mb{R}$ \cite[Theorem 6 b)]{Kruse:Estimates}. However, to apply inductively Moshchevitin's argument in dimension, e.g., $m\times 1$, one would require an estimate of the form
$$\sum_{q=1}^{Q}\frac{1}{q\|q\alpha_{1}\|\dotsc\|q\alpha_{m}\|}\ll_{m}(\log Q)^{m+1}$$
for at least some vectors $(\alpha_{1},\dotsc,\alpha_{m})$. At present, such estimate is only known to hold for multiplicatively badly approximable vectors (to see this, it suffices to apply Abel's summation formula and \cite[Theorem 2.1]{LeVaaler:Sumsof}). Hence, a different method is required.

To prove Proposition \ref{prop:mainprop}, we work directly in a higher-dimensional setting without relying on induction. We generalise a construction introduced by Badziahin and Velani \cite{BadziahinVelani:MultiplicativelyBadly}, in order to produce a multi-dimensional Cantor-like set contained in $\textup{Mad}_{m,n}(C,\bs{\gamma},f,c)$. Such construction requires to count lattice points lying in sets with "hyperbolic spikes". To accomplish this, we use an elementary geometric argument that is the key to the whole proof. The core of this argument can be found in Lemma \ref{lem:hypercount'}. We remark that Badziahin's proof \cite{Badziahin:OnMultiplicatively} of the fact that $\textup{dim}\,\textup{Mad}^{\lambda}(2,1)=1$ for $\lambda>1$ also relies on induction.

Proposition \ref{prop:mainprop} is one $\log$ factor off from the conjecturally optimal result (i.e., the extension of Conjecture \ref{conj:BadzVel} to higher dimension). Specifically, we could not prove that $\textup{Mad}^{\lambda}(m,n)\neq\emptyset$ for $m+n-2\leq \lambda< m+n-1$ in consequence of some overcounting issues araising in the proof of Proposition \ref{prop:mainprop}. Badziahin \cite{Badziahin:OnMultiplicatively} used a rather convoluted strategy to overcome such issues, thus improving on the estimates of Moshchevitin. However, his methods appear quite hard to generalise to a higher-dimensional and/or inhomogeneous setting.

We conclude by saying that it would be equally desirable to prove a dimensional result for the set $\textup{Mad}^{m+n-1}(m,n)$. Unfortunately, the methods used in this paper do not seem powerful enough to obtain such result. Indeed, the (suitably generalised) hypothesis in Badziahin and Velani's \cite[Theorem 4]{BadziahinVelani:MultiplicativelyBadly} does not hold for our construction. An adaptation of \cite[Theorem 4]{BadziahinVelani:MultiplicativelyBadly} to our setting appears equally challenging, due to an obstruction in \cite[Lemma 2]{BadziahinVelani:MultiplicativelyBadly}.

\subsection{Applications}

Let $m,n\in\mb{N}$, let $\bs{Q}\in(0,+\infty)^{n}$, and let $X:=\prod_{j=1}^{n}[-Q_{j},Q_{j}]$. Let also $\bs{L}\in\mb{R}^{m\times n}$ such that the entries $L_{i1},\dotsc,L_{in}\in\mb{R}$ together with $1$ are linearly independent over $\mb{Z}$ for $i=1,\dotsc,m$. We consider the sum
\begin{equation}
S_{\bs{L}}(\bs{Q}):=\sum_{\substack{\bs{q}\in X\cap\mb{Z}^{n}\setminus\{\bs{0}\}}}\prod_{i=1}^{m} \|L_{i}\bs{q}\|^{-1}. \nonumber
\end{equation}
Sums of this shape are of major importance in Diophantine approximation and have extensively been studied (see \cite{Fregoli:OnACounting} for a brief summary). L\^e and Vaaler \cite[Corollary 1.2]{LeVaaler:Sumsof} showed that for $Q:=(Q_{1}\dotsm Q_{n})^{1/n}\geq 1$ it holds
\begin{equation}
S_{\bs{L}}(\bs{Q})\gg_{m,n}Q^{n}(\log Q)^{m}\nonumber
\end{equation}
independently of the choice of the matrix $\bs{L}$. They also asked whether this estimate is sharp, i.e., whether there exist matrices $\bs{L}$ such that
\begin{equation}
S_{\bs{L}}(\bs{Q})\ll_{m,n} Q^{n}(\log Q)^{m}.\nonumber
\end{equation}
In \cite[Theorem 2]{LeVaaler:Sumsof}, they showed that this holds true for multiplicatively badly approximable matrices, but since these matrices are not known to exist, the question remains open. Proposition \ref{prop:mainprop} allows us to find matrices with "relatively small" (even though not optimal) upper bounds.

Let $\phi:[1,+\infty)\to (0,1]$ be a non-increasing function. In \cite[Corollary 1.8]{Fregoli:OnACounting} the author proved that if a matrix $\bs{L}$ is $\phi$-semimultiplicatively badly approximable, i.e., if
\begin{equation}
|\bs{q}|_{\infty}^{n}\prod_{i=1}^{m}\|L_{i}\bs{q}\|\geq\phi(|\bs{q}|_{\infty})\nonumber
\end{equation}
for all $\bs{q}\in\mb{Z}^{n}\setminus\{\bs{0}\}$, then we have
\begin{equation}
\label{eq:condition1'}
\sum_{\substack{\bs{q}\in[-Q,Q]^{n}\\ \cap\ \mb{Z}^{n}\setminus\{\bs{0}\}}}\prod_{i=1}^{m} \|L_{i}\bs{q}\|^{-1}\ll_{m,n}Q^{n}\log\left(\frac{Q^{n}}{\phi(Q)}\right)^{m}+\frac{Q^{n}}{\phi(Q)}\log\left(\frac{Q^{n}}{\phi(Q)}\right)^{m-1}
\end{equation}
for $Q\geq 2$. Since $\prod_{j=1}^{n}\max\{1,|q_{j}|\}\leq|\bs{q}|_{\infty}^{n}$ for all $\bs{q}\in\mb{Z}^{n}$, from Proposition \ref{prop:mainprop} we easily deduce the following.
\begin{cor}
Let $m,n\in\mb{N}$. Then, there exist uncountably many matrices $\bs{L}\in\mb{R}^{m\times n}$ such that
\begin{equation}
\label{eq:sumsof}
S_{\bs{L}}(\bs{Q})\ll_{m,n} Q^{n}(\log Q)^{2m+n-2}
\end{equation}
for all $\bs{Q}=(Q,\dotsc,Q)$ with $Q\geq 2$.
\end{cor}

\noindent Note that the linear independence of the row entries of $\bs{L}$ together with $1$ over $\mb{Z}$ follows from the definition of $\textup{Mad}^{m+n-1}(m,n)$.

This result is not best possible. In particular, by (\ref{eq:condition1'}), we have that for $\bs{L}\in\textup{Mad}^{1+\varepsilon}(2,1)$ it holds
$$S_{\bs{L}}(\bs{Q})\ll_{\varepsilon} Q(\log Q)^{1+\varepsilon}$$
for any $\varepsilon>0$ (such matrices $\bs{L}$ exist since $\textup{dim}\,\textup{Mad}^{\lambda}(2,1)=2$ for $\lambda>1$). Hence, for $m=2,\ n=1$ inequality (\ref{eq:sumsof}) is not sharp. It is also well-known (see \cite{Fregoli:Sumsof}) that set of $1\times n$ matrices $\bs{L}$ such that
$$S_{\bs{L}}(\bs{Q})\ll_{n} Q^{n}\log Q$$
has full Hausdorff dimension in $\mb{R}^{1\times n}$. Thus, (\ref{eq:sumsof}) is again not sharp for $m=1$. However, to the best of our knowledge, for $m\geq 3$ or $m=2,\ n\geq 2$ the existence of matrices satisfying (\ref{eq:sumsof}) was not previously known.
 
\section{Generalised Cantor sets in higher dimension}
\label{sec:Cantor}

In this section we introduce a simple generalisation\footnote{To be precise, our construction is simplified. In \cite{BadziahinVelani:MultiplicativelyBadly}, Badziahin and Velani consider a double-index sequence $\bs{r}=(r_{h,k})$, whereas we consider two one-index sequences $\bs{r}=(r_{k})$ and $\bs{h}=(h_{k})$, since this is enough for our application.} of the one-dimensional construction used by Badziahin and Velani in \cite{BadziahinVelani:MultiplicativelyBadly}. This generalisation will be useful in the proof of Proposition \ref{prop:mainprop}. From now on the word cube will stand for ball in the supremum norm.

Let $l\in\mb{N}$ and let $C$ be a closed cube. For $k\geq 0$ let $\bs{R}:=(R_{k})$ be a sequence of natural numbers, and let $\bs{r}:=(r_{k})$ and $\bs{h}:=(h_{k})$ be sequences of non-negative integers with $0\leq h_{k}\leq k$.

Our goal is to construct a Cantor-like set contained in $C$ depending on the sequences $\bs{R}$, $\bs{r}$, and $\bs{h}$. We denote such set by $\bs{K}(C,\bs{R},\bs{h},\bs{r})$. To this end, we introduce two sequences $\mc{I}_{k}$ and $\mc{J}_{k}$ of cube collections such that each cube in these collections lies in $C$ ($k\geq 0$). We set $\mc{I}_{0}=\mc{J}_{0}:=\{C\}$ and we define $\mc{I}_{k}$ and $\mc{J}_{k}$ by recursion on $k$. We do this in two steps. Suppose that we have constructed $\mc{I}_{h}$ and $\mc{J}_{h}$ for $h=0,\dotsc,k$. Then,\vspace{2mm}
\begin{itemize}
\item[STEP 1] we split each cube $J\in\mc{J}_{k}$ into $R_{k}^{l}$ cubes of equal volume. We call $\mc{I}_{k+1}$ the family of all the cubes obtained via this splitting procedure for $J$ ranging in $\mc{J}_{k}$; note that for $I\in\mc{I}_{k+1}$
$$\textup{edge}(I)=R_{k}^{-l}\textup{edge}(J)\quad\mbox{and}\quad\#\mc{I}_{k+1}=R_{k}^{l}\#\mc{J}_{k};$$\vspace{-0.5mm} 
\item[STEP 2] for each $J\in\mc{J}_{h_{k}}$ we remove from $\mc{I}_{k+1}$ at most $r_{k}$ cubes $I\in\mc{I}_{k+1}$ such that $I\subset J$. We call $\mc{J}_{k+1}$ the family given by the remaining cubes in $\mc{I}_{k+1}$.
\end{itemize}\vspace{2mm}
Finally, we set
\begin{equation}
\mbf{K}(C,\bs{R},\bs{h},\bs{r}):=\bigcap_{k=1}^{\infty}\bigcup_{J\in\mc{J}_{k}}J.\nonumber
\end{equation}
Note that the sequences $\bs{R}$, $\bs{r}$, and $\bs{h}$ do not determine a unique set, but a number of different sets obtained via the procedure described above. This follows from the fact that we did not specify which cubes we remove in the second step (we only gave a bound on their number). We call every set constructed by using the sequences $\bs{R}$, $\bs{r}$, and $\bs{h}$, in the cube $C$, a $(C,\bs{R},\bs{h},\bs{r})$-Cantor set. We also observe that, by construction,
\begin{equation}
\label{eq:Cantor2}
\#\mc{J}_{k+1}\geq R_{k}^{l}\#\mc{J}_{k}-r_{k}\#\mc{J}_{h_{k}}
\end{equation}
for all $k\geq 0$.

Now, the following proposition extends \cite[Theorem 3]{BadziahinVelani:MultiplicativelyBadly}.

\begin{prop}[multidimensional Baziahin-Velani]
\label{prop:Cantor}
Let $\mbf{K}(C,\bs{R},\bs{h},\bs{r})$ be a $(C,\bs{R},\bs{h},\bs{r})$-Cantor set, where $C\subset\mb{R}^{l}$ is a cube, and let
$$t_{k}:=R_{k}^{l}-\frac{r_{k}}{\prod_{i=h_{k}}^{k-1}t_{i}}$$
for $k\geq 1$. If $t_{k}>0$ for all $k\geq 0$, then we have $\mbf{K}(C,\bs{R},\bs{h},\bs{r})\neq\emptyset$.
\end{prop}

The proof is almost straightforward and we give it directly in this section.

\begin{proof}
We shall prove by induction on $k$ that for $k\geq 1$
\begin{equation}
\label{eq:Jk}
\#\mc{J}_{k}\geq t_{k-1}\#\mc{J}_{k-1}.
\end{equation}
The fact that $t_{k}>0$ for all $k$, along with (\ref{eq:Jk}), implies
$$\#\mc{J}_{k}\geq\left(\prod_{h=0}^{k-1}t_{h}\right)\#\mc{J}_{0}>0.$$
Hence, every $(C,\bs{R},\bs{h},\bs{r})$-Cantor set is the intersection of a family of nested compact non-empty sets, and therefore non-empty.

We are left to prove that $\#\mc{J}_{k}\geq t_{k-1}\#\mc{J}_{k-1}$ for all $k\geq 1$. By (\ref{eq:Cantor2}), we have $\#\mc{J}_{1}\geq R_{0}^{l}\#\mc{J}_{0}-r_{0}\#\mc{J}_{0}=t_{0}\#\mc{J}_{0}$, and this proves the case $k=1$.
Now, let us assume that for all $1\leq h\leq k$ it holds $\#\mc{J}_{h}\geq t_{h-1}\#\mc{J}_{h-1}$. Then, in particular, we have
$$\#\mc{J}_{k}\geq\left(\prod_{i=h_{k}}^{k-1}t_{i}\right)\#\mc{J}_{h_{k}}.$$
This, combined with (\ref{eq:Cantor2}), gives
\begin{equation}
 \#\mc{J}_{k+1}\geq R_{k}^{l}\#\mc{J}_{k}-r_{k}\#\mc{J}_{h_{k}}\geq\left(R_{k}^{l}-\frac{r_{k}}{\prod_{i=h_{k}}^{k-1}t_{i}}\right)\#\mc{J}_{k}=t_{k}\#\mc{J}_{k},\nonumber
\end{equation}
whence the claim.
\end{proof}

\section{Proof of Proposition \ref{prop:mainprop}}
\label{sec:Proof}

The strategy is simple enough: by picking suitable parameters, we construct a non-empty $(C,\bs{R},\bs{h},\bs{r})$-Cantor set $\mbf{K}(C,\bs{R},\bs{h},\bs{r})$ lying in $\textup{Mad}_{m,n}\left(C,\bs{\gamma},\log^{*}(x)^{m+n-1},c\right)\setminus\bigcup_{H\in\mathscr{H}}H$.

To do so, we fix a non-decreasing\footnote{We say that a sequence $\{x_{i}\}_{i\in\{0\}\cup\mb{N}}$ of real numbers is non-decreasing if $x_{i}\leq x_{i+1}$ for all $i$.} sequence of integers $\bs{R}=(R_{k})$ with $R_{k}\geq 1$, a sequence of non-negative integers $\bs{h}$, with $0\leq h_{k}\leq k$, and a strictly increasing unbounded function $F:\{0\}\cup\mb{N}\to[1,+\infty)$. In the following technical lemma we specify the values of a sequence $\bs{r}$ (in terms of $c, \ell, F$, $\bs{R}$, and $\bs{h}$) for which there exists a (possibly empty) $(C,\bs{R},\bs{h},\bs{r})$-Cantor set lying in $\textup{Mad}_{m,n}\left(C,\bs{\gamma},\log^{*}(x)^{m+n-1},c\right)\setminus\bigcup_{H\in\mathscr{H}}H$. 

\begin{lem}
\label{lem:MadCantor}
Assume that\vspace{2mm}
\begin{itemize}
\item[$i)$] $2^{m}c<e^{-1}$;\vspace{2mm}
\item[$ii)$] $F(0)= 1$ and $F(k+1)/F(k)\geq e$ for all $k\geq 0$;\vspace{2mm}
\item[$iii)$] $F(k+1)^{2}\log^{*}(F(k+1))^{m+n-1}\leq c\ell^{-1}\prod_{h=0}^{k}R_{h}$ for all $k\geq 0$.\vspace{2mm}
\end{itemize}
Then, there is a $(C,\bs{R},\bs{h},\bs{r})$-Cantor set contained in $\textup{Mad}_{m+n}\left(C,\bs{\gamma},\log^{*}(x)^{m+n-1},c\right)\setminus\bigcup_{H\in\mathscr{H}}H$ with $\bs{r}$ given by
\begin{equation}
\label{eq:rk,hk}
r_{k}:=
\textup{const}(m,n)\left[\mathfrak{f}(c,\ell,\bs{R},\bs{h},k)\prod_{h=h_{k}}^{k}R_{h}^{mn}+\prod_{h=h_{k}}^{k}R_{h}^{mn-1}\right],
\end{equation}
where the factor $\mathfrak{f}(c,\ell,\bs{R},\bs{h},k)$ has the form
\begin{multline}
\label{eq:mathfrakf}
\mathfrak{f}(c,\ell,\bs{R},\bs{h},k):=c\log\left(\frac{1}{2^{m}c}\right)^{m-1}\frac{1}{\log^{*}(F(k))}\log\left(\frac{F(k+1)}{F(k)}\right)^{n-1} \\
\left(\log\left(\frac{F(k+1)}{F(k)}\right)+\ell^{-m}\left(2F(k)^{-m/n}-F(k+1)^{-m/n}\right)\prod_{h=0}^{h_{k}-1}R_{h}^{m}\right),
\end{multline}
and $\textup{const}(m,n)>0$ is a constant only depending on $m$ and $n$. 
\end{lem}

\noindent Lemma \ref{lem:MadCantor} is a key result in our method. Its proof, although quite technical, is essentially based on elementary geometric considerations. We prove Lemma \ref{lem:MadCantor} in Section \ref{sec:proofofMadCantor}.

Now, we need to show that the $(C,\bs{R},\bs{h},\bs{r})$-Cantor set constructed in Lemma \ref{lem:MadCantor} is non-empty. To do so, we use a non-emptiness condition involving the values of the sequence $\bs{r}$.

\begin{lem}
\label{lem:nonempty}
Let $\bs{K}(C,\bs{R},\bs{h},\bs{r})$ be a $(C,\bs{R},\bs{h},\bs{r})$-Cantor set. If for all $k\geq 0$ we have
\begin{equation}
\label{eq:nonempty}
r_{k}\leq\frac{g_{k}}{\max\{2,k\}}\prod_{h=h_{k}}^{k}R_{h}^{mn},
\end{equation}
where $g_{k}:=\max\{2,h_{k}\}/(8\max\{2,k-1\})$, then the set $\bs{K}(C,\bs{R},\bs{h},\bs{r})$ is non-empty.
\end{lem}

\noindent We prove this lemma in Section \ref{sec:proofofnonempty}.

To conclude the proof of Proposition \ref{prop:mainprop}, it is enough to show that both the hypotheses of Lemma \ref{lem:MadCantor} and Lemma \ref{lem:nonempty} simultaneously hold for an appropriate choice of the parameters $c,F,\bs{R},$ and $\bs{h}$. With this in mind, we fix a constant $R>0$, and we set $R_{k}:=R$, $F(k):=R^{k/3}$, and $h_{k}:=\left\lfloor k/(3n) \right\rfloor$ for all $k\geq 0$. Then, we prove that, provided $R$ is large enough, the constant $c$ has enough room to satisfy both the hypotheses of Lemma \ref{lem:MadCantor} and Lemma \ref{lem:nonempty}.

With our choice of $\bs{R},F,$ and $\bs{h}$, condition $ii)$ in Lemma \ref{lem:MadCantor} becomes $R\geq e^{3}$, whereas condition $iii)$ becomes
\begin{equation}
R^{\frac{2(k+1)}{3}}\log^{*}\left(R^{\frac{k+1}{3}}\right)^{m+n-1}\leq c\ell^{-1}R^{k+1},\nonumber
\end{equation}
whence
\begin{equation}
\label{eq:cond3}
\ell R^{-\frac{k+1}{3}}\log^{*}\left(R^{\frac{k+1}{3}}\right)^{m+n-1}\leq c.
\end{equation}
On the other hand, by substituting (\ref{eq:rk,hk}) into (\ref{eq:nonempty}), we obtain
\begin{equation}
\textup{const}(m,n)\left[\mathfrak{f}(c,\ell,\bs{R},\bs{h},k)\prod_{h=h_{k}}^{k}R_{h}^{mn}+\prod_{h=h_{k}}^{k}R_{h}^{mn-1}\right]\leq\frac{g_{k}}{\max\{2,k\}}\prod_{h=h_{k}}^{k}R_{h}^{mn},\nonumber
\end{equation}
which, with our choice of $\bs{R},F,$ and $\bs{h}$, is equivalent to
\begin{equation}
\label{eq:f1}
\mathfrak{f}(c,\ell,R,k)+\frac{1}{R^{\left(k-\left\lfloor\frac{k}{3n}\right\rfloor+1\right)}}\leq\frac{g_{k}\textup{const}(m,n)^{-1}}{\max\{2,k\}}.
\end{equation}
Since $g_{k}$ is bounded away from $0$ for all $k$, by choosing $R$ suitably large in terms of $m$ and $n$, we can ignore the second term at the left-hand side of (\ref{eq:f1}). Hence, we are just left to prove
$$\mathfrak{f}(c,\ell,R,k)\leq\frac{\textup{const}'(m,n)}{\max\{2,k\}},$$
where $\textup{const}'(m,n)$ is a constant only depending on $m$ and $n$. By using (\ref{eq:mathfrakf}), this can be written as
\begin{multline}
\label{eq:cond2}
c\log^{*}\left(\frac{1}{2^{m}c}\right)^{m-1}\frac{1}{\max\{1,k\}}\log^{*}\left(R^{1/3}\right)^{n-1} \\ 
\left(\log^{*}\left(R^{1/3}\right)+\ell^{-m}R^{\frac{-mk}{3n}}\left(2-R^{-\frac{m}{3n}}\right)R^{\left\lfloor\frac{k}{3n}\right\rfloor m}\right)\leq\frac{\textup{const}'(m,n)}{\max\{2,k\}},
\end{multline}
where we ignored a factor of $\log\left(R^{1/3}\right)$ at the denominator, coming from $\log^{*}(F(k))$ for $k\geq 1$. Assuming $\ell<1$, condition (\ref{eq:cond2}) holds if we have
\begin{equation}
\label{eq:cond4}
c\log^{*}\left(\frac{1}{2^{m}c}\right)^{m-1}\leq\textup{const}''(m,n)\ell^{m}\log^{*}\left(R^{1/3}\right)^{-n},
\end{equation}
where $\textup{const}''(m,n)$ is some other positive constant only depending on $m$ and $n$.

To conclude the proof, we pick a small real number $\varepsilon>0$. Since $\log^{*}(1/2^{m}c)^{m-1}\ll_{m,\varepsilon}c^{-\varepsilon}$, condition (\ref{eq:cond4}) is in turn implied by
\begin{equation}
\label{eq:cond5}
c\leq \textup{const}'''(m,n,\varepsilon)\ell^{\frac{m}{1-\varepsilon}}\log^{*}\left(R^{1/3}\right)^{-\frac{n}{1-\varepsilon}},
\end{equation}
where $\textup{const}'''(m,n,\varepsilon)$ is a suitably chosen positive constant only depending on $m$, $n$, and $\varepsilon$. The claim is then proved on noting that (\ref{eq:cond3}) and (\ref{eq:cond5}) can simultaneously hold for a sufficiently large value of $R$.

\section{Proof of Lemma \ref{lem:MadCantor}}
\label{sec:proofofMadCantor}

\subsection{Construction of the Cantor-like set}
For each $P:=(\bs{p},\bs{q})\in\mb{Z}^{m}\times\left(\mb{Z}^{n}\setminus\{\bs{0}\}\right)$, with $gcd(p_{1},\dotsc,p_{m},q_{1},\dotsc,q_{n})=1$, we introduce the following "bad" set.
\begin{multline}
\Delta(P):=\left\{\bs{X}\in\mb{R}^{m\times n}:\prod_{i=1}^{m}\left|X_{i}\bs{q}+\gamma_{i}+p_{i}\right|\leq\frac{c}{\prod(\bs{q})\log^{*}\left(\prod(\bs{q})\right)^{m+n-1}},\right. \\
\left|X_{i}\bs{q}+\gamma_{i}+p_{i}\right|\leq\frac{1}{2}\ i=1,\dotsc,m\Bigg\},
\end{multline}
where we ignore the dependence on $\bs{\gamma}$ and $c$ for simplicity. We also enumerate the hyperplanes in $\mathscr{H}$, indexing them for $k\in\{0\}\cup\mb{N}$. Then, we define the families $\mc{J}_{k}$ of our Cantor-like set so that the intersection of their cubes avoids all the "bad" sets $\Delta(P)$ and the hyperplanes $H_{k}$ for $k\in\mb{N}$. More precisely, for each $J\in\mc{J}_{k}$ we require that $J\cap(\Delta(P)\cup H_{h})=\emptyset$ for all the points $P$ with $\prod(\bs{q})< F(k)$ and all hyperplanes $H_{h}$ with $h\leq k$ (where we assume $H_{0}=\emptyset$). If this condition is satisfied, we have 
\begin{multline}
\label{eq:containment}
\mbf{K}(C,\bs{R},\bs{h},\bs{r})\subset\bigcap_{k=0}^{+\infty}\bigcap_{\prod(\bs{q})<F(k)}\bigcap_{h\leq k}C\setminus(\Delta(P)\cup H_{h}) \\
=C\setminus\left(\bigcup_{P}\Delta(P)\cup\bigcup_{H\in\mathscr{H}}H\right)=\textup{Mad}_{m,n}\left(C,\bs{\gamma},\log^{*}(x)^{m+n-1},c\right)\setminus\bigcup_{H\in\mathscr{H}}H,
\end{multline}
thus showing the claim. Note that (\ref{eq:containment}) holds because the function $F(k)$ is unbounded.

We construct the families $\mc{J}_{k}$ by recursion on $k\geq 0$. For each $k$ we need to ensure
\begin{equation}
\label{eq:P}
J\in\mc{J}_{k}\Rightarrow J\cap\left(\bigcup_{\prod(\bs{q})<F(k)}\Delta(P)\cup\bigcup_{h\leq k}H_{h}\right)=\emptyset.
\end{equation}
If $k=0$, we have $\mc{J}_{0}=\{C\}$ and $F(0)=1$. Therefore, by definition,
$$\bigcup_{\prod(\bs{q})<1}\Delta(P)\cup\bigcup_{h\leq 0}H_{h}=\emptyset.$$ 
This shows that $\mc{J}_{0}$ satisfies (\ref{eq:P}). For $k\geq 1$ we subdivide the points $P\in\mb{Z}^{m}\times\left(\mb{Z}^{n}\setminus\{\bs{0}\}\right)$ into "workable" families. Namely, we define
$$C(k):=\left\{P\in\mb{Z}^{m}\times\left(\mb{Z}^{n}\setminus\{\bs{0}\}\right):F(k-1)\leq \prod(\bs{q})<F(k)\right\}.$$

Suppose that we have constructed the family $\mc{J}_{k}$ in such a way that for $J\in\mc{J}_{k}$ (\ref{eq:P}) holds (note that $\mc{J}_{k}$ can be empty). If $\mc{J}_{k}=\emptyset$, we set $\mc{J}_{k+1}:=\emptyset$, if $\mc{J}_{k}\neq\emptyset$, we proceed as follows. Since any cube in $\mc{I}_{k+1}$ lies within some cube in $\mc{J}_{k}$, it is enough to construct $\mc{J}_{k+1}$ in such a way that if $J\in\mc{J}_{k+1}$ then $J\cap(\Delta(P)\cup H_{k+1})=\emptyset$ for all $P\in C(k+1)$. To define $\mc{J}_{k+1}$, we therefore remove from $\mc{I}_{k+1}$ all the cubes $I$ such that $I\cap(\Delta(P)\cup H_{k+1})\neq\emptyset$  for some $P\in C(k+1)$. This procedure yields a possibly empty Cantor-like set $\bs{K}(C,\bs{R},\bs{h},\bs{r})$, contained in $\textup{Mad}_{m,n}\left(C,\bs{\gamma},\log^{*}(x)^{m+n-1},c\right)\setminus\bigcup_{H\in\mathscr{H}}H$. 

To conclude the proof, we just need to estimate the number of "small" cubes $I\in \mc{I}_{k+1}$ that need to be removed from each "big" cube $J\in\mc{J}_{h_{k}}$ to avoid the sets $\Delta(P)$ for $P\in C(k+1)$ and the hyperplane $H_{k+1}$, and show that such number is smaller than $r_{k}$ defined in (\ref{eq:rk,hk}).

We start by counting the cubes intersecting the sets $\Delta(P)$ for $P\in C(k+1)$. In particular, for a fixed $J\in\mc{J}_{h_{k}}$ it is enough to estimate
$$\#\{I\in\mc{I}_{k+1}:\exists P\in C(k+1)\ I\cap J\cap\Delta(P)\neq\emptyset\}.$$
If $\mc{I}_{k+1}=\emptyset$, there is nothing to prove. Otherwise, we write
\begin{multline}
\{I\in\mc{I}_{k+1}:\exists P\in C(k+1)\ I\cap J\cap\Delta(P)\neq\emptyset\} \\
=\bigcup_{\substack{\bs{q}\in\mb{Z}^{n}\setminus\{\bs{0}\} \\ F(k)\leq\prod(\bs{q})<F(k+1)}}\bigcup_{\substack{P\in C(k+1) \\ \bs{q}(P)=\bs{q}}}\{I\in\mc{I}_{k+1}:I\cap J\cap\Delta(P)\neq\emptyset\},\nonumber
\end{multline}
whence we deduce
\begin{equation}
\label{eq:A+B}
\#\{I\in\mc{I}_{k+1}:\exists P\in C(k+1)\ I\cap J\cap\Delta(P)\neq\emptyset\}\leq\sum_{\substack{\bs{q}\in\mb{Z}^{n}\setminus\{\bs{0}\} \\ F(k)\leq\prod(\bs{q})<F(k+1)}} A(\bs{q})B(\bs{q}),
\end{equation}
where
\begin{equation}
A(\bs{q}):=\max_{\substack{P\in C(k+1) \\ \bs{q}(P)=\bs{q}}}\#\{I\in\mc{I}_{k+1}:I\cap J\cap\Delta(P)\neq\emptyset\}
\end{equation}
and
\begin{equation}
B(\bs{q}):=\#\left\{P\in C(k+1):\bs{q}(P)=\bs{q},\ J\cap\Delta(P)\neq\emptyset\right\}.
\end{equation}
We estimate the factors $A(\bs{q})$ and $B(\bs{q})$ separately.

\subsection{Estimate of $A(\bs{q})$}

To estimate $A(\bs{q})$, we need the following counting result.

\begin{lem}
\label{lem:hypercount'}
Let $\bs{\gamma}'\in\mb{R}^{m}$, $\bs{q}\in\mb{Z}^{n}\setminus\{\bs{0}\}$, and $\varepsilon,T\in(0,+\infty)$, with $\varepsilon/T^{m}<e^{-1}$ (where $e=2.71828\dots$ is the base of the natural logarithm). Let also
$$\mathscr{C}:=\left\{\bs{X}\in\mb{R}^{m\times n}:\prod_{i=1}^{m}|X_{i}\bs{q}+\gamma_{i}'|\leq\varepsilon,\ |X_{i}\bs{q}+\gamma_{i}'|\leq T,\ i=1,\dotsc,m\right\},$$
and let $\mathscr{D}\subset\mb{R}^{m\times n}$ be a cube such that $\mathscr{D}\cap\mathscr{C}\neq\emptyset$. Finally, let $\delta>0$, $\bs{V}\in\mb{R}^{m\times n}$, and $\Lambda$ be the grid $\delta\mb{Z}^{m\times n}+\bs{V}$. Then,
we have
\begin{multline}
\label{eq:hypercount'}
\delta^{mn}\#\{\mbox{tiles } \tau\ \mbox{of the grid } \Lambda: \tau\cap\mathscr{D}\cap\mathscr{C}\neq\emptyset\}\leq 2^{2m-1}\frac{\varepsilon +(T+n|\bs{q}|_{\infty}\delta)^{m}-T^{m}}{|\bs{q}|_{\infty}^{m}} \\
\log^{*}\left(\frac{(T+n|\bs{q}|_{\infty}\delta)^{m}}{\varepsilon +(T+n|\bs{q}|_{\infty}\delta)^{m}-T^{m}}\right)^{m-1}(\textup{edge}(\mathscr{D})+2\delta)^{m(n-1)},
\end{multline}
where a tile is any set of the form $\{\bs{X}\in\mb{R}^{m\times n}:\delta S_{ij}+V_{ij}\leq X_{ij}\leq \delta (S_{ij}+1)+V_{ij},\ i=1,\dotsc,m,\ j=1,\dotsc,n\}$ for some $\bs{S}\in\mb{Z}^{m\times n}$.
\end{lem}

\noindent We prove this result in Section \ref{sec:proofofgemolemma}.

Now, we note that if\vspace{2mm}
\begin{itemize}
\item[$a)$] $\varepsilon\gg_{m,n}(T+n|\bs{q}|_{\infty}\delta)^{m}-T^{m}$;\vspace{2mm}
\item[$b)$] $T\gg_{m,n}|\bs{q}|_{\infty}\delta$;\vspace{2mm}
\item[$c)$] $\textup{edge}(\mathscr{D})\gg_{m,n}\delta$;\vspace{2mm}
\end{itemize}
(\ref{eq:hypercount'}) implies
\begin{equation}
\label{eq:hypercount''}
\#\{\mbox{tiles } \tau\ \mbox{of the grid } \Lambda: \tau\cap\mathscr{D}\cap\mathscr{C}\neq\emptyset\}\ll_{m,n}\frac{\varepsilon}{\delta^{mn}|\bs{q}|_{\infty}^{m}}\log^{*}\left(\frac{T^{m}}{\varepsilon}\right)^{m-1}\textup{edge}(\mathscr{D})^{m(n-1)}.
\end{equation}

This is precisely the assertion that we need to prove the claim. We fix a point $P\in C(k+1)$ and a cube $J\in\mc{J}_{h_{k}}$, and we apply (\ref{eq:hypercount''}) to $\mathscr{C}=\Delta(P)$, $\mathscr{D}=J$, and to the grid $\Lambda$ formed by the cubes $I\in\mc{I}_{k+1}$. We have
$\varepsilon=c\left(\prod(\bs{q})\log^{*}\left(\prod(\bs{q})\right)^{m+n-1}\right)^{-1}$, $T=1/2$, and $\delta=\ell\prod_{h=0}^{k}R_{h}^{-1}$ (note that by hypothesis $\varepsilon/T^{m}<e^{-1}$).

We show that conditions $a)$, $b)$, and $c)$ hold in this specific case. If condition $b)$ is satisfied, then to prove condition $a)$, it is enough to show that $\varepsilon\gg_{m,n}T^{m-1}|\bs{q}|_{\infty}\delta$. By definition of $C(k+1)$ and part $iii)$ in the hypotheses of Lemma \ref{lem:MadCantor}, we have
\begin{equation}
\label{eq:proof6}
\frac{\varepsilon}{|\bs{q}|_{\infty}}=\frac{c}{|\bs{q}|_{\infty}\prod(\bs{q})\log^{*}\left(\prod(\bs{q})\right)^{m+n-1}}\geq\frac{c}{F(k+1)^{2}\log^{*}(F(k+1))^{m+n-1}}\geq \ell\prod_{h=0}^{k}R_{h}^{-1}=\delta.
\end{equation}
Hence, $\varepsilon\gg_{m,n}T^{m-1}|\bs{q}|_{\infty}\delta$, and we have $a)$. Condition $b)$ is equivalent to $1/|\bs{q}|_{\infty}\gg_{m,n}\delta$, which is again implied by (\ref{eq:proof6}).
Finally, condition $c)$ is clearly satisfied since $\textup{edge}(J)\geq\textup{edge}(I)$ for any $I\in\mc{I}_{k+1}$.

Thus, we can apply (\ref{eq:hypercount''}) to obtain
\begin{multline}
\label{eq:A}
A(\bs{q})=\max_{\substack{P\in C(k+1) \\ \bs{q}(P)=\bs{q}}}\#\{I\in\mc{I}_{k+1}:I\cap J\cap\Delta(P)\neq\emptyset\}\ll_{m,n}\frac{c}{|\bs{q}|_{\infty}^{m}\prod(\bs{q})\log^{*}\left(\prod(\bs{q})\right)^{m+n-1}} \\
\log^{*}\left(\frac{\prod(\bs{q})\log^{*}\left(\prod(\bs{q})\right)^{m+n-1}}{2^{m}c}\right)^{m-1}\ell^{-m}\prod_{h=0}^{h_{k}-1}R_{h}^{-m(n-1)}\prod_{h=0}^{k}R_{h}^{mn} \\
\ll_{m,n}\frac{c\log^{*}(1/(2^{m}c))^{m-1}}{|\bs{q}|_{\infty}^{m}\prod(\bs{q})\log^{*}\left(\prod(\bs{q})\right)^{n}}\ell^{-m}\prod_{h=0}^{h_{k}-1}R_{h}^{-m(n-1)}\prod_{h=0}^{k}R_{h}^{mn}.
\end{multline}

\subsection{Estimate of $B(\bs{q})$}

We are now left to estimate $\#\{P\in C(k+1):\bs{q}(P)=\bs{q},\ J\cap\Delta(P)\neq\emptyset\}$ for each given $\bs{q}\in\mb{Z}^{n}\setminus\{\bs{0}\}$ such that $F(k)\leq\prod(\bs{q})<F(k+1)$. For $P\in C(k+1)$ we consider the hyperspace
$$\pi(P):=\{\bs{X}\in\mb{R}^{m\times n}:X_{i}\bs{q}+\gamma_{i}+p_{i}=0,\ i=1,\dotsc,m\},$$
i.e., the "core" of the set $\Delta(P)$. We show that to count the number of points $P\in C(k+1)$ such that $\Delta(P)$ intersects $J$, it is enough to count the number of points $P\in C(k+1)$ such that the thinner set $\pi(P)$ intersects an "inflation" of $J$. In particular, we claim that
\begin{multline}
\label{eq:proof9}
\#\{P\in C(k+1):\bs{q}(P)=\bs{q},\ J\cap\Delta(P)\neq\emptyset\} \\
\leq\#\{P\in C(k+1):\bs{q}(P)=\bs{q},\ J_{\sqrt{m}/|\bs{q}|_{\infty}}\cap\pi(P)\neq\emptyset\},
\end{multline}
where $J_{\sqrt{m}/|\bs{q}|_{\infty}}$ is the "inflation" of the cube $J$ by the quantity $\sqrt{m}/|\bs{q}|_{\infty}$, i.e., the set $\{\bs{X}\in\mb{R}^{m\times n}: \textup{dist}_{\infty}(\bs{X},J)\leq \sqrt{m}/|\bs{q}|_{\infty}\}$. To prove (\ref{eq:proof9}) we show that for any fixed $P$ in the left-hand side of (\ref{eq:proof9}) we have
$$J_{\sqrt{m}/|\bs{q}|_{\infty}}\cap\pi(P)\neq\emptyset.$$
Indeed, for any $\bs{Y}\in\Delta(P)$ we have
$$|Y_{i}\bs{q}+\gamma_{i}+p_{i}|\leq 1/2\ \mbox{for }i=1,\dotsc,m.$$
Moreover, for $i=1,\dotsc,m$ the Euclidean distance in $\mb{R}^{m}$ between the vector $Y_{i}$ and the hyperplane $\{X_{i}\bs{q}+\gamma_{i}+p_{i}=0\}$ is given by $|Y_{i}\bs{q}+\gamma_{i}+p_{i}|/|\bs{q}|_{2}$. Hence, the Euclidean distance in $\mb{R}^{m\times n}$ between the vector $\bs{Y}$ and $\pi(P)$ is at most $\sqrt{m}/(2|\bs{q}|_{2})$. This shows that for any point $\bs{Y}\in\Delta(P)$, we have
$$\textup{dist}_{2}(\bs{Y},\pi(P))\leq\frac{\sqrt{m}}{2|\bs{q}|_{2}}.$$
Since $J\cap\Delta(P)\neq\emptyset$, we deduce
$$\textup{dist}_{\infty}(J,\pi(P))\leq\textup{dist}_{2}(J,\pi(P))\leq\textup{dist}_{2}(J\cap\Delta(P),\pi(P))\leq\frac{\sqrt{m}}{2|\bs{q}|_{2}}\leq\frac{\sqrt{m}}{2|\bs{q}|_{\infty}}.$$
Hence, by definition of distance, $J_{\sqrt{m}/|\bs{q}|_{\infty}}\cap\pi(P)\neq\emptyset$, whence the claim.

We are now left to bound the right-hand side in (\ref{eq:proof9}). To do this, we note that the distance between two hyperspaces $\pi(P)$ and $\pi(P')$ with same $\bs{q}$ is at least $1/(n|\bs{q}|_{\infty})$. Indeed, assume that $X_{i},X_{i}'\in\mb{R}^{n}$ satisfy $X_{i}\bs{q}+\gamma_{i}+p_{i}=0$ and $X_{i}'\bs{q}+\gamma_{i}+p_{i}'=0$, with $p_{i}\neq p_{i}'$. Then, by the Cauchy-Schwartz inequality, we have
\begin{equation}
\textup{dist}_{\infty}(X_{i},X_{i}')\geq\frac{\textup{dist}_{2}(X_{i},X_{i}')}{\sqrt{n}}\geq\frac{|(X_{i}-X_{i}')\bs{q}|}{\sqrt{n}|\bs{q}|_{2}}\geq\frac{|p_{i}-p_{i}'|}{\sqrt{n}|\bs{q}|_{2}}\geq \frac{1}{n|\bs{q}|_{\infty}}.\nonumber
\end{equation}
This shows that for all $\bs{X}\in\pi(P)$ and $\bs{X}'\in\pi(P')$ we have
\begin{equation}
\label{eq:planedist}
\textup{dist}_{\infty}(\bs{X},\bs{X}')\geq\frac{1}{n|\bs{q}|_{\infty}}.
\end{equation}

Now, if $J_{\sqrt{m}/|\bs{q}|_{\infty}}\cap\pi(P)\neq\emptyset$ for some $P$, by a dimensional argument\footnote{Observe that $\pi(P)$ must intersect the boundary of $J$, hence some $(mn-1)$-dimensional face $F$ of $J$. The intersection of $\pi(P)$ with the hyperspace generated by $F$ has dimension at least $\textup{dim}(F)+\textup{dim}(\pi(P))-\textup{dim}(F+\pi(P))\geq \textup{dim}(\pi(P))-1$. Hence, we have a hyperspace of dimension $\pi(P)-1$ intersecting a cube (F) of dimension $mn-1$. The argument can be run inductively.}, the hyperspace $\pi(P)$ must intersect at least one $m$-dimensional face of the cube $J_{\sqrt{m}/|\bs{q}|_{\infty}}$. For each $P$ such that $J_{\sqrt{m}/|\bs{q}|_{\infty}}\cap\pi(P)\neq\emptyset$ we select a point $Q(P)$ on an $m$-dimensional face of $J_{\sqrt{m}/|\bs{q}|_{\infty}}$ lying in $\pi(P)$. We know, by (\ref{eq:planedist}), that all such points are at least at a distance of $1/(n|\bs{q}|_{\infty})$ away from each other in the supremum distance. To evaluate their number, we fix any $m$-dimensional face $E$ of $J_{\sqrt{m}/|\bs{q}|_{\infty}}$ and we enlarge it by $1/(2n|\bs{q}|_{\infty})$ in all directions, i.e., we consider the set $E_{1/(2n|\bs{q}|_{\infty})}$. Then, for each intersection point $Q(P)\in\pi(P)\cap J_{\sqrt{m}/|\bs{q}|_{\infty}}$ we take an $mn$-dimensional cube of edge $1/(n|\bs{q}|_{\infty})$ centred at $Q(P)$. All these cubes are contained in the "inflated" face $E_{1/(2n|\bs{q}|_{\infty})}$ and they all have disjoint interiors. Comparing the volume of these cubes and the volume of the inflated face, we find
\begin{multline}
\label{eq:newedge}
\#\{P\in C(k+1):E\cap\pi(P)\neq\emptyset\}\left(\frac{1}{n|\bs{q}|_{\infty}}\right)^{mn} \\
\leq\textup{Vol}\left(E_{1/(2n|\bs{q}|_{\infty})}\right)=\left(\textup{edge}\left(J_{\sqrt{m}/|\bs{q}|_{\infty}}\right)+\frac{1}{n|\bs{q}|_{\infty}}\right)^{m}\left(\frac{1}{n|\bs{q}|_{\infty}}\right)^{m(n-1)} \\
=\left(\textup{edge}\left(J\right)+\frac{1+2n\sqrt{m}}{n|\bs{q}|_{\infty}}\right)^{m}\left(\frac{1}{n|\bs{q}|_{\infty}}\right)^{m(n-1)}.
\end{multline}
Since the number of $m$-dimensional faces of a cube only depends on $m$, from (\ref{eq:proof9}) and (\ref{eq:newedge}) we deduce
\begin{multline}
\label{eq:B}
B(\bs{q})=\#\{P\in C(k+1):\bs{q}(P)=\bs{q},\ J\cap\Delta(P)\neq\emptyset\}\leq \#\{P\in C(k+1): \\ \bs{q}(P)=\bs{q},\ J_{\sqrt{m}/(2|\bs{q}|_{\infty})}\cap\pi(P)\neq\emptyset\}
\ll_{m,n}\left(|\bs{q}|_{\infty}\textup{edge}\left(J\right)+1\right)^{m}\ll_{m,n}|\bs{q}|_{\infty}^{m}\ell^{m}\prod_{h=1}^{h_{k}-1}R_{h}^{-m}+1.
\end{multline}

\subsection{Conclusion}

To conclude the proof of Lemma \ref{lem:MadCantor}, we combine (\ref{eq:A+B}), (\ref{eq:A}), and (\ref{eq:B}) to obtain
\begin{multline}
\#\{I\in\mc{I}_{k+1}:\exists P\in C(k+1)\ I\cap J\cap\Delta(P)\neq\emptyset\} \\
\ll_{m,n}\sum_{\substack{\bs{q}\in\mb{Z}^{n}\setminus\{\bs{0}\} \\ F(k)\leq\prod(\bs{q})<F(k+1)}}\left(\frac{c\log^{*}(1/(2^{m}c))^{m-1}}{|\bs{q}|_{\infty}^{m}\prod(\bs{q})\log^{*}\left(\prod(\bs{q})\right)^{n}}\ell^{-m}\prod_{h=0}^{h_{k}-1}R_{h}^{-m(n-1)}\prod_{h=0}^{k}R_{h}^{mn}\right) \\
\left(|\bs{q}|_{\infty}^{m}\ell^{m}\prod_{h=1}^{h_{k}-1}R_{h}^{-m}+1\right).\nonumber
\end{multline}
Hence, by using the fact that $|\bs{q}|_{\infty}^{m}\geq\prod(\bs{q})^{m/n}$, we find
\begin{multline}
\label{eq:34}
\#\{I\in\mc{I}_{k+1}:\exists P\in C(k+1)\ I\cap J\cap\Delta(P)\neq\emptyset\}\ll_{m,n} c\log^{*}\left(\frac{1}{2^{m}c}\right)^{m-1}\prod_{h=h_{k}}^{k}R_{h}^{mn} \\
\frac{1}{\log^{*}\left(F(k)\right)^{n}}\sum_{\substack{\bs{q}\in\mb{Z}^{n}\setminus\{\bs{0}\} \\ F(k)\leq\prod(\bs{q})<F(k+1)}}\frac{1}{\prod(\bs{q})}\left(1+\frac{\ell^{-m}}{\prod(\bs{q})^{m/n}}\prod_{h=0}^{h_{k}-1}R_{h}^{m}\right).
\end{multline}
Now, a simple integration shows that
$$\sum_{\substack{\bs{q}\in\mb{Z}^{n}\setminus\{\bs{0}\} \\ F(k)\leq\prod(\bs{q})<F(k+1)}}\prod(\bs{q})^{-1}\ll_{n}\log^{*}(F(k+1))^{n-1}\log^{*}\left(\frac{F(k+1)}{F(k)}\right),$$
and
$$\sum_{\substack{\bs{q}\in\mb{Z}^{n}\setminus\{\bs{0}\} \\ F(k)\leq\prod(\bs{q})<F(k+1)}}\prod(\bs{q})^{-1-m/n}\ll_{n}\log^{*}(F(k+1))^{n-1}\left(2F(k)^{-m/n}-F(k+1)^{-m/n}\right).$$
Therefore, from (\ref{eq:34}) and from $\log^{*}(F(k+1))^{n-1}\leq\log^{*}(F(k+1)/F(k))^{n-1}\log^{*}(F(k))^{n-1}$, we deduce
\begin{multline}
\label{eq:)}
\#\{I\in\mc{I}_{k+1}:\exists P\in C(k+1)\ I\cap J\cap\Delta(P)\neq\emptyset\}\ll_{m,n} \\
c\log^{*}\left(\frac{1}{2^{m}c}\right)^{m-1}\frac{1}{\log^{*}\left(F(k)\right)}\log^{*}\left(\frac{F(k+1)}{F(k)}\right)^{n-1} \\
\left(\log^{*}\left(\frac{F(k+1)}{F(k)}\right)+\ell^{-m}\left(2F(k)^{-m/n}-F(k+1)^{-m/n}\right)\prod_{h=0}^{h_{k}-1}R_{h}^{m}\right)\prod_{h=h_{k}}^{k}R_{h}^{mn}.
\end{multline}

We are now left to count all the cubes in $\mc{I}_{k+1}$ lying in $J$ and intersecting the hyperplane $H_{k+1}$. From the set of cubes $I\in\mc{I}_{k+1}$ such that $I\cap J\cap H_{k+1}\neq\emptyset$, we select a maximal subset $S$ of pairwise disjoint cubes (with disjoint boundary). For each of these cubes $I$ we pick a point lying in $I\cap H_{k+1}$. The points that we picked are, by construction, at least $\textup{edge}(I)$ distant from each other in the supremum norm. Then, we take $(mn-1)$-dimensional cubes in $H_{k+1}$ of edge $\textup{edge}(I)$ around each such point. By construction, these cubes are disjoint. Comparing the volume of the union of the cubes with the volume of the set $(J\cap H_{k+1})$ inflated in the Euclidean distance by the quantity $\textup{diam}(I)$ in $H_{k+1}$, i.e., the set $\left\{\bs{X}\in H_{k+1}:\textup{dist}_{2}(\bs{X},J\cap H_{k+1})\leq\textup{diam}(I)\right\}$, we find
$$\# S\cdot\textup{edge}(I)^{mn-1}\ll_{m,n}\left(\textup{diam}(J\cap H_{k+1})+\textup{diam}(I)\right)^{mn-1}\ll_{m,n}\textup{edge}(J)^{mn-1},$$  
whence
\begin{equation}
\label{eq:f2}
\#\left\{I\in\mc{I}_{k+1}:I\cap J\cap H_{k+1}\neq\emptyset\right\}\ll_{m,n}\# S\ll_{m,n}\left(\frac{\textup{edge}(J)}{\textup{edge}(I)}\right)^{mn-1}=\prod_{h=h_{k}}^{k}R_{h}^{mn-1}.
\end{equation}
Combining (\ref{eq:)}) and (\ref{eq:f2}), the proof of Lemma \ref{lem:MadCantor} is concluded.

\section{Proof of Lemma \ref{lem:nonempty}}
\label{sec:proofofnonempty}

We show by induction on $k$ that
\begin{equation}
\label{eq:proof1}
t_{k}\geq R_{k}^{mn}\left(1-\frac{1}{\max\{2,k\}}\right)>0
\end{equation}
for all $k\geq 0$. By Proposition \ref{prop:Cantor}, this is enough to prove the claim. If $k=0$ we have
\begin{equation}
t_{0}=R_{0}^{mn}-r_{0}\geq R_{0}^{mn}-\frac{g_{0}}{2}R_{0}^{mn}.\nonumber
\end{equation}
Hence, the base case is proved, given that $g_{0}=1/8$. Now, assume that (\ref{eq:proof1}) holds for $0\leq h\leq k$. Then, we have
\begin{align}
\label{eq:g1}
t_{k+1} & =R_{k+1}^{mn}-\frac{r_{k+1}}{\prod_{i=h_{k+1}}^{k}t_{i}}\nonumber \\
 & \geq R_{k+1}^{mn}-\frac{r_{k+1}}{\prod_{i=h_{k+1}}^{k}R_{h}^{mn}\left(1-\max\{2,i\}^{-1}\right)}.
\end{align}
Moreover, for $k\geq 0$
$$\prod_{i=h_{k+1}}^{k}\left(1-\max\{2,i\}^{-1}\right)\geq\frac{\max\left\{2,h_{k+1}\right\}-1}{4\max\{2,k\}}\geq\frac{\max\{2,h_{k+1}\}}{8\max\{2,k\}}=g_{k+1}.$$
Hence, by (\ref{eq:g1}) and by the hypothesis, we deduce
\begin{equation}
t_{k+1}\geq R_{k+1}^{mn}-\frac{g_{k+1}^{-1} r_{k+1}}{\prod_{h=h_{k+1}}^{k}R_{h}^{mn}}\geq R_{k+1}^{mn}\left(1-\frac{1}{\max\{2,k+1\}}\right).\nonumber
\end{equation}

\section{Proof of Lemma \ref{lem:hypercount'}}
\label{sec:proofofgemolemma}

For a set $\mathscr{A}\subset\mb{R}^{m\times n}$ we denote by $\mathscr{A}_{\delta}$ the "inflation" of $\mathscr{A}$ by the quantity $\delta$, i.e., the set $\{\bs{X}\in\mb{R}^{m\times n}:\textup{dist}_{\infty}(\bs{X},\mathscr{A})\leq\delta\}$. First, we show that
\begin{equation}
\delta^{mn}\#\{\mbox{tiles } \tau\ \mbox{of the lattice } \Lambda: \tau\cap\mathscr{D}\cap\mathscr{C}\neq\emptyset\}\leq\textup{Vol}(\mathscr{D}_{\delta}\cap\mathscr{C}_{\delta}).
\end{equation}
This follows from the fact that for any tile $\tau$ of $\Lambda$ we have
\begin{equation}
\label{eq:geomlemma1}
\tau\cap\mathscr{D}\cap\mathscr{C}\neq\emptyset\Rightarrow\tau\subset\mathscr{D}_{\delta}\cap\mathscr{C}_{\delta}.
\end{equation}
To see why (\ref{eq:geomlemma1}) holds, it is enough to observe that for all points $P\in\tau\cap\mathscr{C}$ and all points $Q\in\tau$ we have $\textup{dist}_{\infty}(P,Q)\leq\delta$. Hence, $\tau\subset\mathscr{C}_{\delta}$. The same is true for $\mathscr{D}$, whence (\ref{eq:geomlemma1}).

To conclude the proof, we need to estimate $\textup{Vol}(\mathscr{D}_{\delta}\cap\mathscr{C}_{\delta})$. By definition, if $\bs{X}\in\mathscr{C}_{\delta}$, then there is some $\bs{X}'\in\mathscr{C}$ such that $\textup{dist}_{\infty}(\bs{X},\bs{X}')\leq\delta$. Hence,
\begin{multline}
\label{eq:geomlemma1.5}
\prod_{i=1}^{m}|X_{i}\bs{q}+\gamma_{i}'|\leq\prod_{i=1}^{m}\left(|X'_{i}\bs{q}+\gamma_{i}'|+n|\bs{q}|_{\infty}\delta\right) \\
\leq \prod_{i=1}^{m}|X'_{i}\bs{q}+\gamma_{i}'|+\sum_{I\subsetneq\{1,\dotsc,m\}}\left(\prod_{i\in I}|X'_{i}\bs{q}+\gamma_{i}'|\right)(n|\bs{q}|_{\infty}\delta)^{m-|I|}\leq\varepsilon+(T+n|\bs{q}|_{\infty}\delta)^{m}-T^{m}.
\end{multline}
Now, let $\bs{\mu}$ be the centre of the cube $\mathscr{D}$. Then, by (\ref{eq:geomlemma1.5}), we have
\begin{multline}
\label{eq:geomlemma2}
\mathscr{D}_{\delta}\cap\mathscr{C}_{\delta}\subset \\
\left\{\bs{X}\in\mb{R}^{m\times n}:\begin{cases}
\prod_{i=1}^{m}|X_{i}\bs{q}+\gamma_{i}'|\leq\varepsilon +(T+n|\bs{q}|_{\infty}\delta)^{m}-T^{m} & \\
|X_{i}\bs{q}+\gamma_{i}'|\leq T+n|\bs{q}|_{\infty}\delta & i=1,\dotsc,m \\
|X_{ij}-\mu_{ij}|\leq\textup{edge}(\mathscr{D})/2+\delta & i=1,\dotsc,m,\ j=1,\dotsc,n
\end{cases}\right\}.
\end{multline}

Without loss of generality, we can assume $|\bs{q}|_{\infty}=|q_{n}|$. We proceed by considering the linear transformation $\xi:\mb{R}^{m\times n}\to\mb{R}^{m\times n}$ defined by
$$\xi(\bs{X})_{ij}=\begin{cases}
X_{ij} & \mbox{if }j\neq n \\
X_{i}\bs{q} & \mbox{if }j=n
\end{cases}.$$ 
Under the action of $\xi$, the right-hand side of (\ref{eq:geomlemma2}) is sent into a subset of
\begin{equation}
\label{eq:geomlemma3}
\left\{\bs{X}\in\mb{R}^{m\times n}:\begin{cases}
\prod_{i=1}^{m}|X_{in}+\gamma_{i}'|\leq\varepsilon +(T+n|\bs{q}|_{\infty}\delta)^{m}-T^{m} & \\
|X_{in}+\gamma_{i}'|\leq T+n|\bs{q}|_{\infty}\delta & i=1,\dotsc,m \\
|X_{ij}-\mu_{ij}|\leq\textup{edge}(\mathscr{D})/2+\delta & i=1,\dotsc,m\ j=1,\dotsc,n-1
\end{cases}\right\}.
\end{equation}
Now, the determinant of $\xi$ is $|q_{n}|^{m}=|\bs{q}|_{\infty}^{m}\neq 0$, so $\xi$ is a bijective linear transformation. Therefore, to obtain an estimate of $\textup{Vol}(\mathscr{D}_{\delta}\cap\mathscr{C}_{\delta})$, it is enough to estimate the volume of the set in (\ref{eq:geomlemma3}). A simple integration shows that the volume of this set is bounded from above by
$$2^{2m-1}(\varepsilon +(T+n|\bs{q}|_{\infty}\delta)^{m}-T^{m})\log^{*}\left(\frac{(T+n|\bs{q}|_{\infty}\delta)^{m}}{\varepsilon +(T+n|\bs{q}|_{\infty}\delta)^{m}-T^{m}}\right)^{m-1}(\textup{edge}(\mathscr{D})+2\delta)^{m(n-1)}.$$

Hence, $\textup{Vol}(\mathscr{D}_{\delta}\cap\mathscr{C}_{\delta})$ is bounded from above by this quantity divided by the absolute value of the determinant of $\xi$, that is $|q_{n}|^{m}=|\bs{q}|_{\infty}^{m}.$

\addcontentsline{toc}{section}{\bibname}
\bibliographystyle{plain}
\bibliography{Bibliography}

\begin{thebibliography}{10}

\bibitem{Badziahin:OnMultiplicatively}
D.~Badziahin.
\newblock On multiplicatively badly approximable numbers.
\newblock {\em Mathematika}, 59(No. 1):31--55, 2013.

\bibitem{BadziahinVelani:MultiplicativelyBadly}
D.~Badziahin and S.~Velani.
\newblock Multiplicatively badly approximable numbers and generalised {C}antor
  sets.
\newblock {\em Advances in {M}athematics}, 228(No. 5):2766--2796, 2011.

\bibitem{Beresnevich:Classic}
V.~Beresnevich and S.~Velani.
\newblock Classical metric {D}iophantine approximation revisited: the
  {K}hintchine-{G}roshev {T}heorem.
\newblock {\em {I}nt. {M}ath. {R}es. {N}ot.}, 2010(No. 1):69--86.

\bibitem{Bugeaud:Multiplicative}
Y.~Bugeaud.
\newblock `{M}ultiplicative {D}iophantine approximation', {D}ynamical systems
  and {D}iophantine approximation.
\newblock {\em Proc. Conf. Inst. H. Poincar\'{e} (Soci\'{e}t\'{e}
  Math\'{e}matique de France, Paris)}, pages 105--125, 2009.

\bibitem{BugeaudMoshchevitin:BadlyApp}
Y.~Bugeaud and N.~Moshchevitin.
\newblock Badly approximable numbers and {L}ittlewood-type problems.
\newblock {\em Math. Proc. Cambridge Phil. Soc.}, 150(No. 2):215--226, 2011.

\bibitem{Fregoli:Sumsof}
R.~Fregoli.
\newblock Sums of reciprocals of fractional parts.
\newblock {\em {I}nt. {J}. {N}umber {T}heory}, 15(No. 4):789--797, 2019.

\bibitem{Fregoli:OnACounting}
R.~Fregoli.
\newblock On a counting theorem for weakly admissible lattices.
\newblock {\em Int. Mat. Res. Not.}, rnaa102, 2020.

\bibitem{Gallagher:MetricSimultaneous}
P.~Gallagher.
\newblock Metric simultaneous {D}iophantine approximation.
\newblock {\em J. London Math. Soc.}, 37(No. 1):387--390, 1962.

\bibitem{Kruse:Estimates}
A.~H. Kruse.
\newblock Estimates of {$\sum_{k=1}^{N}k^{-s}\langle kx\rangle^{-t}$}.
\newblock {\em Trans. Amer. Math. Soc.}, 110:493--518, 1964.

\bibitem{LeVaaler:Sumsof}
T.~H. L\^{e} and J.~D. Vaaler.
\newblock Sums of products of fractional parts.
\newblock {\em Proc. London Math. Soc.}, 111(No. 3):561--590, 2014.

\bibitem{Moshchevitin:BadlyApp}
N.~Moshchevitin.
\newblock Badly approximable numbers related to the {L}ittlewood conjecture.
\newblock {\em arXiv:0810.0777 [math.NT]}, 2008.

\bibitem{Peck:Simultaneous}
L.~G. Peck.
\newblock Simultaneous rational approximations to algebraic numbers.
\newblock {\em Bull. A.M.S.}, 67:197--201, 1961.

\bibitem{PollingtonVelani:OnaProblem}
A.~Pollington and S.~Velani.
\newblock On a problem in simultaneous {D}iophantine approximation:
  {L}ittlewood's conjecture.
\newblock {\em Acta Math.}, 66:29--40, 2000.

\bibitem{Schmidt:BadlyApprox}
W.~M. Schmidt.
\newblock Badly approximable systems of linear forms.
\newblock {\em J. Number Theory}, 1:139--154, 1969.

\bibitem{Sprindzuk:Metric}
V.~G. Sprind\v{z}uk.
\newblock {\em Metric theory of {D}iophantine approximations (in {R}ussian)}.
\newblock Nauka, 1977.

\end{thebibliography}

\end{document}